\documentclass[11pt]{amsart}

\usepackage{amsmath,amssymb,amscd,enumitem}
\usepackage{amsrefs}
\usepackage{xcolor}

\usepackage{amsmath,amssymb}
\usepackage{float}

\usepackage{amsxtra, amsmath}
\usepackage{amssymb, amscd}
\usepackage{graphicx}
\usepackage{mathrsfs}

\usepackage{mathtools}

\newcommand{\CC}{\mathbb{C}}

\newtheorem{thm}{Theorem}[section]
\newtheorem{prop}[thm]{Proposition}
\newtheorem{lemma}[thm]{Lemma}
\newtheorem{cor}[thm]{Corollary}

\theoremstyle{definition}
\newtheorem{definition}[thm]{Definition}

\theoremstyle{remark}
\newtheorem{remark}[thm]{Remark}

\usepackage{dsfont}
\allowdisplaybreaks

\title 
{Singular solutions of the matrix Bochner problem: the $N$-dimensional cases }

\author{Ignacio Bono Parisi}
\author{Ines Pacharoni}

\address{CIEM-FaMAF, Universidad Nacional de C\'ordoba, C\'ordoba~5000, Argentina}
\email{ines.pacharoni@unc.edu.ar}
\email{ignacio.bono@unc.edu.ar}

\thanks{This paper was partially supported by SeCyT-UNC, CONICET, PIP 11220200102031CO}

\begin{document}
\begin{abstract}
In the theory of matrix-valued orthogonal polynomials, there exists a longstanding problem known as the Matrix Bochner Problem: the classification of all \( N \times N \) weight matrices \( W(x) \) such that the associated orthogonal polynomials are eigenfunctions of a second-order differential operator. In \cite{CY18}, Casper and Yakimov made an important breakthrough in this area, proving that, under certain hypotheses, every solution to this problem can be obtained as a bispectral Darboux transformation of a direct sum of classical scalar weights.

In the present paper, we construct three families of weight matrices \( W(x) \) of size \( N \times N \), associated with Hermite, Laguerre, and Jacobi weights, which can be considered `singular' solutions to the Matrix Bochner Problem because they cannot be obtained as a Darboux transformation of classical scalar weights.

\end{abstract} 
\maketitle

\section{Introduction}

Orthogonal matrix polynomials are sequences of matrix-valued polynomials that are pairwise orthogonal with respect to a matrix-valued inner product defined by an $N\times N$ weight matrix $W(x)$.
The theory of these matrix-valued orthogonal polynomials goes back to  M. G. Krein in \cite{K49} and \cite{K71}. In \cite{D97}, R. Dur\'an  initiated the study of matrix-valued orthogonal polynomials that furthermore,  are eigenfunctions of certain second-order symmetric differential operators.

In the scalar case, corresponding to $N=1$, a seminal result from S. Bochner \cite{B29} classifies all sequences of orthogonal polynomials on the real line that are eigenfunctions of a second-order linear differential operator. According to this result, after an affine transformation, the only weights with these properties are the classical weights $e^{-x^{2}}$, $e^{-x}x^{\alpha}$, and $(1-x)^\alpha(1 + x)^\beta$ of Hermite, Laguerre, and Jacobi, respectively. 

The analog problem of finding weight matrices \( W(x) \) of size \( N \times N \) such that the associated sequence of orthogonal matrix polynomials are eigenfunctions of a second-order matrix differential operator is often referred to as the {Matrix Bochner Problem} (MBP from now on). 
The first nontrivial solutions to the MBP, for \( N = 2 \), were provided by different methods, by A. Durán and F. A. Grünbaum, I. Pacharoni, and J. Tirao, in \cite{DG04}, \cite{G03}, \cite{GPT01}, and \cite{GPT02}. Over the past twenty years, many additional examples have been discovered, contributing significantly to this subject. See, for instance, \cite{GPT05}, \cite{CG06}, \cite{DG07}, \cite{PT07}, \cite{P08}, \cite{PR08}, \cite{PZ16}, \cite{KPR12}, and \cite{BP24}.

Given a weight $W$, the MBP turns out to be  equivalent to the existence of a second-order differential operator in the algebra $\mathcal{D}(W)$, of  differential operators  with polynomial coefficients  such that 
$$P_n \cdot D=\Lambda_n(D) P_n, \quad \text{ with } \Lambda_{n}(D) \in \operatorname{Mat}_{N}(\mathbb{C}), \, \text{ for all } n\in\mathbb{N}_0.$$
Here $\{P_n\}_{n\in \mathbb{N}_0}$ denotes a sequence of matrix-valued orthogonal polynomials with respect to the weight $W$.

 In \cite{CY18}, R. Casper and M. Yakimov using techniques from noncommutative algebra and representation theory for the algebra $\mathcal{D}(W)$ formulated a general framework for studying the MBP,  which resulted in a significant breakthrough in this area. 
 
 They proved that solutions to MBP, under the additional condition that $\mathcal{D}(W)$ is a full algebra,  arise from  Darboux transformations of a direct sum of classical Hermite, Laguerre, or Jacobi weights, \cite{CY18}*{Theorem 1.3}.

While this result represents a substantial advancement, it gives rise to some remaining open questions: Are there weights that are solutions to the Bochner problem but cannot be obtained as Darboux transformations of direct sums of classical weights? Alternatively, are there weight matrices whose algebra of differential operators is not full? In this case, do these weights exist in all dimensions?

The main goal of this paper is to construct  
three new families of $N\times N$ irreducible weight matrices that are solutions to the MBP and cannot be obtained as Darboux transformations of classical scalar weights. 
 Our families of weights $W(x)$ are associated with each of the classical weight families of  Hermite, Laguerre, or Jacobi. Many of the results and the proofs in the paper apply to all families. Therefore, we have chosen to present our three families together, separating proofs and explicit expressions when necessary. 

\smallskip

We now present the unified construction of the weight matrices 
$W(x)$ for the Hermite, Laguerre, and Jacobi cases. These matrices are built using classical scalar weights and a nilpotent matrix $A$.
\begin{definition}\label{W-definition}
Let  $\widetilde W(x)= w_{1}(x) \oplus \cdots \oplus w_{N}(x)$ be a direct sum of classical scalar weights of the same type (Hermite, Laguerre, or Jacobi), such that $w_{i}(x)w_{j}(x)^{-1}$ {\em is not a rational function} for all $i\neq j$. 
We consider the weight matrix given by 
  \begin{equation*} 
      W(x) =  T(x) \widetilde W(x) T(x)^*,
  \end{equation*} 
  where 
 \begin{equation*} 
 T(x) = I+ Ax\, ,   \qquad  \text { } \quad A = \sum_{j=1}^{[N/2]}a_{2j-1}E_{2j-1,2j} + \sum_{j=1}^{[(N-1)/2]}a_{2j}\,E_{2j+1,2j},  
 \end{equation*} 
and $ a_{j}\in \mathbb{R}$, $a_j\neq 0$. Explicitly, if $w_b(x)= e^{-x^2+2bx}$, $w_{\alpha}(x)=e^{-x}x^{\alpha}$, and $w_{\alpha, \beta}(x)=(1-x)^\alpha(1 + x)^\beta$, are the classical scalar weights of Hermite, Laguerre, and  Jacobi then we have
\begin{enumerate}
    \item {\em For Hermite-type:}  \label{WH}
    $$  W_H(x) =  T(x) \widetilde W_H(x) T(x)^* \, , \qquad
    \widetilde{W}_{H}(x) = w_{b_{1}}(x) \oplus \cdots \oplus w_{b_{N}}(x),  \quad  \text{ with } b_{j} \neq b_{i} .
   $$
  \item {\em For Laguerre-type:}   \label{WL}
  $$  W_L(x) =  T(x) \widetilde W_L(x) T(x)^*\, , \qquad  \widetilde{W}_{L}(x) = w_{\alpha_{1}}(x) \oplus \cdots \oplus w_{\alpha_{N}}(x), 
   \quad  \text{ with }  \alpha_{j} - \alpha_{i} \notin \mathbb{Z}.  
    $$
    \item {\em For Jacobi-type:} \label{WJ}     
    $$   W_J(x) =  T(x) \widetilde W_J(x) T(x)^*\, ,\quad   \widetilde{W}_{J}(x) = w_{\alpha_1,\beta_1}(x) \oplus \cdots \oplus w_{\alpha_N,\beta_N}(x),$$ 
with $\alpha_i-\alpha_j\not \in \mathbb Z$ or $\beta_i-\beta_j\not \in \mathbb Z$, and 
    $\alpha_1+\beta_1 = \alpha_j+ \beta_j +1+(-1)^j $, for all $j=1,\dots, N$. 
\end{enumerate}
\end{definition}

\begin{remark}
    The condition  $\alpha_1+\beta_1 = \alpha_j+ \beta_j +1+(-1)^j$ in \eqref{WJ} is necessary to ensure that the weight $W_J$ admits a second-order differential operator in the algebra $\mathcal D(W_J)$.
\end{remark}

\begin{remark}
    The matrix  $A=\displaystyle \sum_{j=1}^{[N/2]}a_{2j-1}E_{2j-1,2j} + \sum_{j=1}^{[(N-1)/2]}a_{2j}\,E_{2j+1,2j}$, is a 2-step nilpotent matrix. Thus, the exponential of $A$ 
  is given by $e^{Ax} =I+Ax=  T(x)$. 
\end{remark}

\smallskip

Having defined the new families of weight matrices, we now present the main results of this paper. 

\begin{thm}\label{main-thm} 
If \( W = W_H \), \( W_L \), or \( W_J \) are the weight matrices introduced in Definition \ref{W-definition}, then:
\begin{enumerate}
    \item[i)] \( W \) is a solution to the Matrix Bochner Problem.
    \item[ii)] \( W \) is not a Darboux transformation of a direct sum of classical scalar weights.
\end{enumerate}
\end{thm}

In the case when \( N = 2 \), the Hermite weights \( W_H \) were previously considered in our earlier work \cite{BP23}. 
We proved that the algebra \( \mathcal{D}(W) \) is a polynomial algebra in a second-order differential operator \( D \). Consequently, it cannot be a full algebra, implying that it is not a Darboux transformation of a direct sum of classical weights. The present work significantly generalizes these results, extending them to arbitrary sizes and to the full range of Hermite, Laguerre, and Jacobi weight families.

\smallskip

The starting point in the proof of the results in Theorem \ref{main-thm} is to study the right Fourier algebra $\mathcal{F}_R(W)$ associated to the weight function $W$. It is a larger algebra than $\mathcal{D}(W)$, that can be described as the algebra of differential operators with polynomial coefficients, which possess an adjoint with respect to $W$, and this adjoint is also a differential operator with polynomial coefficients.
See equations \eqref{Fourier def} and \eqref{fourier algebra} for the formal definition.  

By analyzing the Fourier algebras of the weights  $\widetilde W(x)= w_{1}(x) \oplus \cdots \oplus w_{N}(x)$ and $W(x)=T(x) \widetilde W(x) T(x)^*$, under the hypothesis   that $w_i(x)w_j(x)^{-1}$ is not a rational function,
we  establish (Theorem \ref{algebras}) 
 \begin{equation*} \label{relationD(W)}
     \mathcal{D}(W) \subseteq T(x)\mathcal{D}(\widetilde{W})T^{-1}(x).
 \end{equation*}
 
We show that $\mathcal D(\widetilde W)$ is an algebra of diagonal operators with entries that are polynomials in the classical  Hermite, Laguerre, or Jacobi differential operator. In particular, this implies  that 
$\mathcal D(W)$ is a commutative algebra (see Proposition \ref{weight}). 

The key point in proving that the weights \(W = W_H, W_L,\) and \(W_J\) cannot be obtained through a Darboux transformation of a direct sum of classical scalar weights (Theorem \ref{W-noDarboux}) 
is that any differential operator in \(\mathcal{D}(W)\) is of even order and has a leading coefficient \(F_{2m}\) that is a multiple of \(\rho(x)^m I\). Here, \(\rho(x)\) equals \(1\), \(x\), or \(1-x^2\), depending on whether \(W\) corresponds to the Hermite, Laguerre, or Jacobi weights, respectively.

If $\delta$ is the classical scalar second-order differential operator for the weight \(w\), then
the operator \(\widetilde{D} = \operatorname{diag}(\delta_{1}, \ldots, \delta_{N})\) belongs to \(\mathcal{D}(\widetilde W)\).
We prove that  there  exists  a constant matrix \(\widetilde{K}\) such that the second-order differential operator \(D = T(x)\,\big(\widetilde{D} + \widetilde{K}\big)\,T(x)^{-1}\) belongs to \(\mathcal{D}(W)\) (see Proposition \ref{solutionDgeneral}).
Furthermore, we obtain the explicit expressions of the second-order differential operators in the algebra $\mathcal D(W)$ (Theorem \ref{boch sol}).
Finally, after confirming the existence of a second-order differential operator $D \in \mathcal{D}(W)$ we establish that every differential operator in $\mathcal{D}(W)$ can be expressed as a polynomial in $D$ (see Theorem \ref{poly alg}).

\smallskip 

The families of weight matrices introduced in Theorem \ref{main-thm}   can be considered {\it singular solutions} to the MBP because they do not conform to the initial characterization proposed in \cite{CY18}. This challenges the prevailing framework and exposes some limitations of the results in \cite{CY18}. Rather than contradicting the earlier findings, our examples serve as a complementary extension, enriching the understanding of the MBP by demonstrating the existence of solutions outside the previously understood scope.

\section{Preliminaries}\label{Backgr}
We begin with a brief review of the basic theory of orthogonal matrix polynomials, the Darboux transformations, and Bochner's problem.

\subsection{Classical orthogonal polynomials}
We know that the only scalar weights supported on an interval in $\mathbb{R}$ for which their associated orthogonal polynomials are eigenfunctions of some second-order differential operator are, up to an affine change of coordinates, 
the classical weights $e^{-x^{2}}$, $e^{-x}x^{\alpha}$, and $(1-x)^\alpha(1 + x)^\beta$ of Hermite, Laguerre, and Jacobi, respectively. We refer to all of them as the {\em classical scalar weights}.  
In all cases, this second-order differential operator can be written as 
 \begin{equation*}
     \delta= \partial ^2 \rho(x) + \partial \, \frac{ (w(x) \rho(x))' }{w(x)},
 \end{equation*} 
where $\delta$ and $\rho$ are as given in Table \ref{weights-table}. Also, it is a known fact that $\mathcal D(w)$, the algebra of all differential operators that have a sequence of these classical orthogonal polynomials as eigenfunctions, is a polynomial algebra in the differential operator $\delta$.  See \cite{M05}.

\begin{remark}
Throughout this paper, we will consider the classical Hermite weights, together with the affine transformation $w_{b}(x) = e^{-x^{2}+2bx}$ that will be necessary for the construction of our weight matrices.
\end{remark}

\begin{table}[H] 
\begin{center} 
\small
    \begin{tabular}{|c|c|c|c|c|}
     \hline    Type   & weight   & support & $\rho(x)$  &  
     differential operator 
     \\    \hline &&&& \\  Hermite    & $w_b(x)=e^{-x^2+2bx}$ & $(-\infty, \infty)$ & 1 &   $\delta =  \partial^{2} -2 \partial (x-b)$ \\ \hline & & & & \\   Laguerre     & $w_\alpha(x)= e^{-x} x^\alpha$ & $ (0,\infty )$  & $x$  &  $\delta =  \partial^{2}x  + \partial(\alpha+1-x)$ \\ \hline  & & & & \\
     Jacobi   & $w_{\alpha,\beta}(x)= (1-x)^\alpha (1+x)^\beta $ & $(-1,1)$ & $1-x^2$ & $\delta = \partial^{2} (1-x^2)  + \partial (\beta-\alpha-x(\alpha+\beta+2))$ \\ & & & & \\\hline
    \end{tabular}\caption{Classical scalar weights. } \label{weights-table}
\end{center} 
\end{table}

\subsection{Matrix-valued orthogonal polynomials and  the algebra $\mathcal D(W)$}\label{sec-MOP}

 Let $W=W(x)$ be a weight matrix of size $N$ on the real line, that is, a complex $N\times N$ matrix-valued smooth function on a (possibly unbounded) interval $(x_0,x_1)$ such that $W(x)$ is positive definite almost everywhere and with finite moments of all orders. Let $\operatorname{Mat}_N(\mathbb{C})$ be the algebra of $N\times N$ complex matrices and let $\operatorname{Mat}_N(\mathbb{C}[x])$ be the algebra of polynomials in the indeterminate $x$ with coefficients in $\operatorname{Mat}_N(\mathbb{C})$. We consider the following Hermitian sesquilinear form in the linear space $\operatorname{Mat}_N(\mathbb{C}[x])$
\begin{equation*}
  \langle P, Q \rangle =  \langle P, Q \rangle_W = \int_{x_0}^{x_1} P(x) W(x) Q(x)^*\, dx.
\end{equation*}

Given a weight matrix $W$ one can construct sequences 
$\{Q_n\}_{n\in\mathbb{N}_0}$ of matrix-valued orthogonal polynomials. That is, each $Q_n$ is a polynomial of degree $n$ with a nonsingular leading coefficient, and $\langle Q_n, Q_m \rangle = 0$ for $n \neq m$.
We observe that there is a unique sequence of monic orthogonal polynomials $\{P_n\}_{n\in\mathbb{N}_0}$ in $\operatorname{Mat}_N(\mathbb{C}[x])$.
By following a standard argument (see \cite{K49} or \cite{K71}) one shows that the monic orthogonal polynomials $\{P_n\}_{n\in\mathbb{N}_0}$ satisfy a three-term recursion relation
\begin{equation}\label{ttrr}
    x P_n(x)=P_{n+1}(x) + B_{n}P_{n}(x)+ C_nP_{n-1}(x), \qquad n\in\mathbb{N}_0,
\end{equation}
where $P_{-1}=0$ and $B_n, C_n$ are matrices depending on $n$ and not on $x$.

Throughout this paper, we consider differential operators of the form $D = \sum_{j=0}^{m} \partial^j F_{j}(x)$, where $\partial = \frac{d}{dx}$ and $F_{j}(x)$ is a matrix-valued function. These operators act {\em on the right-hand side} of a matrix-valued function as follows
$$P(x) \cdot D = \sum_{j=0}^{m} \partial^j(P)(x) F_{j}(x).$$

\noindent 
We consider the algebra of all differential operators with polynomial coefficients 
$$\operatorname{Mat}_{N}(\Omega[x])=\Big\{D = \sum_{j=0}^{m} \partial^{j}F_{j}(x) \, : F_{j} \in \operatorname{Mat}_{N}(\mathbb{C}[x]) \Big \}.$$
\noindent 
More generally, when necessary, we will also consider $\operatorname{Mat}_{N}(\Omega[[x]])$, the set of all differential operators with coefficients in $\mathbb{C}[[x]]$, the ring of power series with complex coefficients.

\begin{prop}[\cite{GT07}, Propositions 2.6 and 2.7]\label{eigenvalue-prop}
  Let $W=W(x)$ be a weight matrix of size $N\times  N$ and let $\{P_n\}_{n\geq 0}$ be  the sequence of monic orthogonal polynomials in $\operatorname{Mat}_N(\mathbb{C}[x])$. If $D$ is a differential operator of order $s$ such that
  $$P_n\cdot D=\Lambda_n P_n, \qquad \text{for all } n\in\mathbb{N}_0,$$
  with $\Lambda_n\in \operatorname{Mat}_N(\mathbb{C})$, then
 $F_i(x)=\sum_{j=0}^i x^j F_j^i$, $F_j^i \in \operatorname{Mat}_N(\mathbb{C})$, is a polynomial and $\deg(F_i)\leq i$. Moreover, $D$ is determined by the sequence $\{\Lambda_n\}_{n\geq 0}$ and
 \begin{equation}\label{eigenvaluemonicos}
   \Lambda_n=\sum_{i=0}^s [n]_i F_i^i, \qquad \text{for all } n\geq 0,
 \end{equation}
    where $[n]_i=n(n-1)\cdots (n-i+1)$, and  $[n]_0=1$.
\end{prop}
 
Given a weight matrix $W$ and  $\{P_n\}_{n\in \mathbb{N}_0}$  a sequence of matrix-valued orthogonal polynomials with respect to $W$, we consider the algebra 
\begin{equation}\label{algDW}
  \mathcal D(W)=\{D\in \operatorname{Mat}_{N}(\Omega[x])\, : \, P_n \cdot D=\Lambda_n(D) P_n, \, \text{ with } \Lambda_n(D)\in \operatorname{Mat}_N(\mathbb{C}), \text{ for all }n\in\mathbb{N}_0\}.
\end{equation}
We observe that the definition of $\mathcal D(W)$ depends only on the weight matrix $W$ and not on the particular sequence of orthogonal polynomials. (See \cite{GT07}, Corollary 2.5).

\begin{prop} [\cite{GT07}, Proposition 2.8]\label{prop2.8-GT}
For each $n\in\mathbb{N}_0$, the mapping $D\mapsto \Lambda_n(D) $ is a representation of $\mathcal D(W)$ in $\operatorname{Mat}_N(\mathbb{C})$. 
Moreover, the sequence of representations $\{\Lambda_n\}_{n\in\mathbb{N}_0}$ separates the elements of $\mathcal D(W)$, i.e.
if $\Lambda_n(D_1)=\Lambda_n(D_2)$ for all $n\geq 0$ then $D_1=D_2$. 
\end{prop}

The {\em formal adjoint} on $\operatorname{Mat}_{N}(\Omega([[x]])$, denoted  by $\mbox{}^*$, is the unique involution extending Hermitian conjugation on $\operatorname{Mat}_{N}(\mathbb C[x])$ and mapping $\partial I$ to $-\partial I$. 
The {\em formal $W$-adjoint} of $ \mathfrak{D}\in \operatorname{Mat}_{N}(\Omega[x])$ is the differential operator $\mathfrak{D}^{\dagger} \in \operatorname{Mat}_{N}(\Omega[[x]])$ defined by
$$\mathfrak{D}^{\dagger}:= W(x)\mathfrak{D}^{\ast}W(x)^{-1},$$
where $\mathfrak{D}^{\ast}$ is the formal adjoint of $\mathfrak D$. 
An operator $\mathfrak{D}\in \operatorname{Mat}_{N}(\Omega[x])$ is called {\em $W$-adjointable} if there exists $\widetilde {\mathfrak{D}} \in \operatorname{Mat}_{N}(\Omega[x])$, such that
$$\langle P\cdot \mathfrak{D},Q\rangle=\langle P,Q\cdot \widetilde{\mathfrak{D}}\rangle,$$
for all $P,Q\in \operatorname{Mat}_N(\mathbb{C}[x])$. Then we say that the operator $\widetilde{\mathfrak D}$ is the $W$-adjoint of $\mathfrak D $.

\

In this paper, we work with weight matrices $W(x)$ that are 'good enough'. Specifically, these weights are such that for each integer $n \geq 0$, the $n$-th derivative $W^{(n)}(x)$ decreases exponentially at infinity, and there exists a scalar polynomial $p_n(x)$ such that $W^{(n)}(x)p_n(x)$ has finite moments. See \cite{CY18}*{Section 2.2}. With this in mind, the following proposition holds.

\begin{prop}[\cite{CY18}, Prop. 2.23] \label{adjuntas}  
If $\mathfrak{D} \in \operatorname{Mat}_{N}(\Omega[x])$ is $W$-adjointable, then $\mathfrak{D}^{\dagger}$ is the $W$-adjoint of $\mathfrak D$, i.e.
 $$\langle \, P\cdot \mathfrak{D},Q\, \rangle=\langle \, P,\, Q\cdot {\mathfrak{D}}^\dagger \rangle,$$ for all $P,Q\in \operatorname{Mat}_N(\mathbb{C}[x])$.
 \end{prop}
 
For $\mathfrak D= \sum_{j=0}^n \partial ^j F_j \in \operatorname{Mat}_{N}(\Omega[x]) $,  the formal $W$-adjoint of $\mathfrak D$ is given by $\mathfrak{D}^\dagger= \sum_{k=0}^n \partial ^k G_k$, with
\begin{equation}\label{daga}
    G_k=  \sum_{j = 0}^{n-k}(-1)^{n-j} \binom{n-j}{k}
    (WF_{n-j}^*)^{(n-k-j)} W^{-1} , \qquad \text{for } 0\leq k\leq n.  
 \end{equation}

\ 

We say that a differential operator $D \in \operatorname{Mat}{N}(\Omega[x])$ is $W$-{\em symmetric} if $\langle P \cdot D, Q \rangle = \langle P, Q \cdot D \rangle$ for all $P, Q \in \operatorname{Mat}{N}(\mathbb{C}[x])$. In other words, $D$ is $W$-adjointable with $D = D^{\dagger}$.

We have  that   the set $\mathcal S(W)$, of all $W$-symmetric operators in $\mathcal D(W)$, satisfies
$$\mathcal D(W)= \mathcal S (W)\oplus i \mathcal S (W).$$

Furthermore, if $D = \sum_{j=0}^{m} \partial^{j} F_{j}$ is a $W$-symmetric operator with $F_{j}$ being a matrix-valued polynomial of degree less than or equal to $j$, then $D$ belongs to the algebra $\mathcal{D}(W)$ (see \cite{GT07}, Proposition 2.10).

\smallskip

\subsection{Bispectral Darboux transformations}
We recall the notion of right Fourier algebra associated with a weight matrix $W(x)$ as given in \cite{CY18}. 
We consider the space of all semi-infinite sequences of matrix-valued rational functions 
$$ \mathcal{P} =\{ P:\mathbb{C}\times \mathbb{N}_0 \longrightarrow M_N(\mathbb{C}) \, : \, P(x,n)  \text{ is a rational function of $x$, for each fixed $n$}  \}.$$ 
We denote by $\operatorname{Mat}_{N}(\mathcal{S})$ the algebra of all  operators of the form 
 $\mathscr{M}= \sum_{j=-\ell}^{k} A_{j}(n)\, \delta^j$, $A_{j}(n)\in\operatorname{Mat}_{N}(\mathbb{C})$, with its left action defined as
\begin{equation}\label{discreteop}(\mathscr{M}\cdot P)(x,n)= \sum_{j=-\ell}^{k} A_{j}(n)(\delta^j \cdot P)(x,n)= \sum_{j=-\ell}^k A_j(n) P(x,n+j).
\end{equation}

Let  $P_n(x)=P(x,n)$ be the 
sequence of monic orthogonal polynomials for a weight matrix $W$. The right Fourier algebra associated with $W$  is defined  by 
\begin{equation} \label{Fourier def}
\begin{split}
    \mathcal F_R(W)= \mathcal F_R(P) &=\{ \mathfrak D\in \operatorname{Mat}_{N}(\Omega[x])\, : \, \exists \, \mathscr{M} \in \operatorname{Mat}_{N}(\mathcal{S})\text{ such that } P(x,n)\cdot \mathfrak D=\mathscr{M}\cdot P(x,n) \}.
\end{split}
\end{equation}

\smallskip

From \cite{CY18}, Theorem 3.7, we recall the  following explicit description of the right Fourier algebra associated with a weight matrix $W$ 
\begin{equation} \label{fourier algebra}
\mathcal{F}_{R}(W) = \{ \mathfrak{D}\in \operatorname{Mat}_{N}(\Omega[x]):\mathfrak{D} \text{ is } W\text{-adjointable and } \mathfrak{D}^{\dagger}\in \operatorname{Mat}_{N}(\Omega[x]) \}.
\end{equation}

\begin{definition}\label{Darbou-transf-def}
Let $W(x)$ and $\widetilde{W}(x)$ be weight matrices, and let $P_n(x)$ and $\widetilde P_n(x)$ be the respective associated sequences of monic orthogonal polynomials. We say that the sequence $\widetilde{P}_n(x)$ is a bispectral Darboux transformation of $P_n(x)$ if there exist differential operators $\mathfrak{D},\widetilde{\mathfrak{D}} \in \mathcal{F}_{R}(W)$, polynomials $F(x),\widetilde{F}(x)$, and sequences of matrices $C_n, \widetilde{C}_n$ which are nonsingular for almost every $n$ and satisfy
$$C_n\widetilde{P}_n(x) = P_n(x) \cdot \mathfrak{D}F(x)^{-1} \text{ and } \quad  \widetilde{C}_n P_n(x) = \widetilde{P}_n(x) \cdot \widetilde{F}(x)^{-1}\widetilde{\mathfrak{D}}.$$
We say that $\widetilde{W}(x)$ is a {\em
Darboux transformation } of $W(x)$ if $\widetilde{P}_n(x)$ is a bispectral Darboux transformation of $P_n(x)$. 
\end{definition}  

\section{ Fundamental properties of the new weight matrices $W$} \label{Sec newfamilies}

In this section, we establish the fundamental properties of the new weight matrices 
$W$ introduced earlier in Definition \ref{W-definition}. These properties will be essential in proving that these weights cannot be obtained through Darboux transformations of classical weights.

To begin, we focus on the right Fourier algebra of diagonal weight matrices $\widetilde{W}$, where the entries are scalar weights satisfying the condition that $w_i(x)w_j(x)^{-1}$ are not rational functions.

\begin{prop} \label{diag}
    Let $\widetilde{W}(x) = w_{1}(x) \oplus \cdots \oplus w_{N}(x)$ be a direct sum of scalar weights supported on the same interval such that $w_{i}(x)w_{j}(x)^{-1}$ is not a rational function for all $i \not = j$. Then the right Fourier algebra of $\widetilde W$ is a diagonal algebra, i.e., 
    $$\mathcal{F}_{R}(\widetilde{W}) \subseteq \left \{ \sum_{j=0}^{n} \partial^{j}F_{j}(x) \, : \, F_{j}(x) = \operatorname{diag}(p_{1}^{j}(x),\ldots, p_{N}^{j}(x)), \text{ with } p_{i}^{j}\in \mathbb{C}[x] \right \}.$$
\end{prop}
\begin{proof}

We consider $\mathcal{A} = \sum_{j=0}^{n}\partial^{j}F_{j}(x) \in \mathcal{F}_{R}(\widetilde{W})$. The coefficients $F_{n-j}(x) = \sum_{1\leq s,t \leq N} p_{s,t}^{j}(x)E_{s,t}$ are polynomial functions, where $E_{s,t}$ denote the matrices with a $1$ in the $(s,t)$-entry and $0$ elsewhere. We observe that $\mathcal{A}$ is $\widetilde{W}$-adjointable with
$$\mathcal{A}^{\dagger} = \widetilde{W}(x)\mathcal{A}^{\ast}\widetilde{W}(x)^{-1} = \sum_{j=0}^{n}\partial^{j}G_{j}(x)\in \operatorname{Mat}_{N}(\Omega[x]).$$
Using \eqref{daga}, with $k=n$, we find that
$$G_{n}(x) = (-1)^{n}\widetilde{W}(x)F_{n}(x)^{\ast}\widetilde{W}(x)^{-1} = (-1)^{n} \sum_{1\leq s,t \leq N}w_{s}(x)w_{t}(x)^{-1}\overline{p^0_{t,s}}(x)E_{s,t}.$$
Since $G_{n}$ is a matrix-valued polynomial and $w_s(x)w_t(x)^{-1}$ is not a rational function, we conclude that $p^{0}_{s,t} = 0$ for all $s \neq t$, implying that the non-diagonal entries of $F_n$ are also zero.

Now we proceed by induction, by assuming that $F_{n-j}$ is a diagonal matrix for all $0 \leq j \leq k$, and proving that $F_{n-k-1}$ is also a diagonal matrix. In fact, from \eqref{daga}, we have that
\begin{equation*}
    \begin{split}
        G_{n-k-1} & = \sum_{j = 0}^{k}(-1)^{n-j} \binom{n-j}{n-k-1}
        (\widetilde{W}F_{n-j}^*)^{(k+1-j)} \widetilde{W}^{-1} +(-1)^{n-k-1}\widetilde{W}F_{n-k-1}^{\ast}\widetilde{W}^{-1} \\
        & = \sum_{j = 0}^{k}(-1)^{n-j} \binom{n-j}{n-k-1} \sum_{s=1}^{N}(w_{s}\overline{p_{s, s}^{j}})^{(k+1-j)}w_{s}^{-1}E_{s,s} + (-1)^{n-k-1}\sum_{1\leq s,t \leq N}w_{s}w_{t}^{-1}\overline{p^{k+1}_{t,s}}E_{s,t}. 
    \end{split}
\end{equation*}
Therefore, $p_{s,t}^{k+1} = 0$ for all $s \neq t$ and thus $F_{n-k-1}$ is a diagonal matrix. This concludes the proof.
\end{proof}

Proposition \ref{diag} shows that the right Fourier algebra \(\mathcal{F}_{R}(\widetilde{W})\) consists of diagonal matrices. This result allows us to fully characterize the algebra \(\mathcal{D}(\widetilde{W})\) as diagonal operators, where each entry corresponds to a polynomial in the classical second-order differential operator of the respective scalar weight.

\begin{cor}\label{alg diag}
    Let $\widetilde{W}(x) = w_{1}(x) \oplus \cdots \oplus w_{N}(x)$ be a direct sum of classical scalar weights supported on the same interval, such that $w_{i}(x)w_{j}(x)^{-1}$ is not a rational function for all $i \not = j$.  Then 
    \begin{equation*}
     \mathcal{D}(\widetilde{W}) = \left \{ \operatorname{diag}(r_{1}(\delta_{1}),\ldots,r_{N}(\delta_{N})) \, \biggm| \, r_{j} \in \mathbb{C}[x], \text{ and } \delta_{j} = \partial^{2}\rho(x)  + \partial \frac {(w_{j}(x) \rho(x))'}{w_{j}(x)} \right \}. 
 \end{equation*}
 In particular, we have that 
 $$\mathcal{D}(\widetilde{W}) = \mathcal{D}(w_{1}) \oplus \cdots \oplus \mathcal{D}(w_{N}) = \mathbb{C}[\delta_{1}] \oplus \cdots \oplus \mathbb{C}[\delta_{N}].$$
\end{cor} 
\begin{proof}
    The statement follows from the inclusion $\mathcal{D}(\widetilde{W}) \subseteq \mathcal{F}_{R}(\widetilde{W})$ and  the observation that the algebra $\mathcal{D}(w_{j})$ is a polynomial algebra on the differential operator $\delta_{j}$.
\end{proof}

\begin{remark}\label{leadingDWtilde}  If 
    $\mathcal{A}$ is a  differential  operator   in  $\mathcal{D}(\widetilde{W})$  of order $2m$,  then its leading coefficient is of the form $$F_{2m}= \rho(x)^{m} \operatorname{diag} (c_{1},\ldots,c_{N}),$$  for some $c_{1},\ldots,c_{N} \in \mathbb{C}$. 
 \end{remark}

\smallskip

Our families of weight matrices, defined in Definition \ref{W-definition}, are of the form 
\begin{equation*}
W(x) = T(x) \widetilde{W}(x) T(x)^{\ast},
\end{equation*}
where $\widetilde{W}= \widetilde{W}_H, \widetilde{W}_L$ or $\widetilde{W}_J$, and $T$ is an exponential matrix given by $T(x) = I + Ax$, with $A$ a two-step nilpotent matrix  defined by 
$$ A = \sum_{j=1}^{[N/2]}a_{2j-1}E_{2j-1,2j} + \sum_{j=1}^{[(N-1)/2]}a_{2j}\,E_{2j+1,2j}.$$

In the following proposition, we will explore the relationship between the right Fourier algebras $\mathcal{F}_{R}(W)$ and $\mathcal{F}_{R}(\widetilde{W})$. We will use weight matrices related through $W(x) = T(x)\widetilde{W}(x)T(x)^{\ast}$, where $T(x)$ and $T(x)^{-1}$ are matrix-valued polynomials, as in the case of our weights where $T(x) = e^{Ax}=I+xA$.

\begin{prop} \label{fourier}
Let $T(x)$ and $T(x)^{-1}$ be matrix-valued polynomials functions and let 
 $W$ and $\widetilde{W}$ be weight matrices with $W(x) =T(x)\widetilde{W}(x)T(x)^{\ast}$. 
 Then the right Fourier algebras satisfy 
 $$\mathcal{F}_{R}(W) = T(x)\mathcal{F}_{R}(\widetilde{W})T^{-1}(x).$$
 Moreover, an operator $\mathcal{A} \in \mathcal{F}_{R}(W)$ is $W$-symmetric if and only if $\mathcal B= T(x)^{-1}\mathcal{A}T(x) \in \mathcal{F}_{R}(\widetilde{W})$ is a $\widetilde{W}$-symmetric operator.
\end{prop}
\begin{proof}
  Let $\mathcal{A} \in \mathcal{F}_{R}(W)$ and 
    $\mathcal{B} = T^{-1}(x)\mathcal{A}T(x)$. 
    Firstly, we observe that $\mathcal B $ has polynomial coefficients because  $T(x)$ and $T^{-1}(x)$ are polynomial functions and $\mathcal{A}\in  \operatorname{Mat}_{N}(\Omega[x])$. 

The formal $\widetilde W$-adjoint of $\mathcal B$ is \begin{equation}\label{Bdaga}
    \mathcal{B}^{\dagger}  
= \widetilde W(x)\mathcal{B}^{\ast} \widetilde W^{-1}(x) = T^{-1}(x)W(x)\mathcal{A}^* W(x)^{-1}T(x)
= T^{-1}(x)\mathcal{A}^{\dagger}T(x), 
\end{equation}
where $\mathcal{A}^{\dagger}$ denotes the $W$-adjoint of $\mathcal A$. 
  In particular $\mathcal{B}^\dagger \in \operatorname{Mat}_{N}(\Omega[x])$. 
    To show that 
    $\mathcal{B}$ is $\widetilde{W}$-adjointable, 
    let $P,Q \in \operatorname{Mat}_{N}(\mathbb{C}[x])$. We have that
    \begin{align*}
            \langle P\cdot \mathcal{B}, Q \rangle_{\widetilde{W}} &
            = \int_{x_{0}}^{x_{1}}(P(x)T(x)^{-1})\cdot \mathcal{A}W(x)(Q(x)T(x)^{-1})^{\ast}dx  = \langle (PT^{-1})\cdot \mathcal{A} , QT^{-1} \rangle_{W},
 \intertext{ and since $\mathcal{A}$ is 
    $W$-adjointable, we get}
             &= \langle PT^{-1}, (QT^{-1})\cdot \mathcal{A}^{\dagger} \rangle_{W} 
            = \int_{x_{0}}^{x_{1}} P(x) \widetilde{W}(x)(Q(x)\cdot \big(T(x)^{-1}\mathcal{A}^{\dagger}T(x)\big )^{\ast} dx \\
            & = \langle P, Q\cdot (T^{-1}\mathcal{A}^{\dagger}T)\rangle_{\widetilde{W}}.
     \end{align*}
    This proves that $\mathcal{B}$ is $\widetilde{W}$-adjointable and  $\mathcal{B} \in \mathcal{F}_{R}(\widetilde{W})$. Therefore, $\mathcal{F}_{R}(W) \subseteq T(x)\mathcal{F}_{R}(\widetilde{W})T^{-1}(x).$ 
Similarly, we get that $\mathcal{F}_{R}(\widetilde W) \subseteq T^{-1}(x)\mathcal{F}_{R}(W)T(x)$. 

To establish the final assertion in the proposition, we recall that an operator $\mathcal{A} \in \mathcal{F}_{R}(W)$ is $W$-symmetric if and only if $\mathcal{A} = \mathcal{A}^\dagger$. By using \eqref{Bdaga}, we observe that
$$\mathcal{B}^{\dagger} = T^{-1}(x)\mathcal{A}^{\dagger}T(x) = T^{-1}(x)\mathcal{A}T(x) = \mathcal{B},$$
thus concluding the proof.
\end{proof}

\smallskip
The matrix $A$, involved in the definition of our weights $W$, belongs to the following linear subspace of $\operatorname{Mat}_N(\mathbb C)$:

\begin{equation}\label{A-vectorspace}
    \mathcal A = \left \{ \sum_{j=1}^{[N/2]}{m_j} E_{2j-1,2j}+ \sum_{j=1}^{[(N-1)/2]}{n_j} E_{2j+1,2j} \, : m_j, n_j\in \CC  \right \}.
\end{equation}
  Observe that if $A_1, A_2\in \mathcal A$ and  $L$ is a diagonal matrix then 
\begin{equation*}\label{M space}
    A_1A_2=0\, ,  \qquad LA_1 \in \mathcal A \, ,  \qquad \text{ and } \quad  A_1L\in \mathcal A.
\end{equation*}

\smallskip

There is a deep relationship between the algebras $\mathcal D(W)$ and $\mathcal D(\widetilde W)$, for the weights $W$ and $\widetilde W$.

\begin{thm}\label{algebras}
Let $\widetilde{W}$ be a weight matrix with diagonal right Fourier algebra, and define $W(x) = T(x)\widetilde{W}(x)T^{*}(x)$, where  
 $T(x) = e^{Ax} $ and  $A\in \mathcal A$.
 Then
 $$\mathcal{D}(W) \subset T(x)\mathcal{D}(\widetilde{W})T^{-1}(x).$$
\end{thm}

\begin{proof} Since $\mathcal{D}(W) = \mathcal{S}(W) \oplus i\mathcal{S}(W)$, it is enough to prove that $\mathcal{S}(W) \subset T(x)\mathcal{S}(\widetilde{W})T^{-1}(x).$
 Let $D  = \sum_{j=0}^{n}\partial^{j}G_{j} \in \mathcal{D}(W)$ be a $W$-symmetric operator of order $n$. From Proposition \ref{fourier}, we have that $D = T(x)\widetilde{D}T^{-1}(x)$ for some $\widetilde{W}$-symmetric operator $\widetilde{D} = \sum_{j=0}^{n}\partial^{j}F_{j} \in \mathcal{F}_{R}(\widetilde{W})$, the right Fourier algebra of $\widetilde{W}$. 
    
    Since $ D\in \mathcal D(W) $, the coefficients $G_{j}$ are polynomials matrices of degree less than or equal to $j$. They satisfy  
\begin{equation}\label{coeff}
\begin{split}
G_{n}(x) & = F_{n}(x)+x [A,F_{n}(x)], \\
G_j(x) & = F_j(x)+ x[A,F_j(x)]+(j+1)A F_{j+1}(x), \qquad \text{for $0\leq j \leq n-1$},
\end{split}
\end{equation}
since $F_j$ are diagonal matrices and $AF_jA=0$.

Thus, the (diagonal) matrix $F_{j}(x)$ coincides  with the main diagonal of $G_{j}(x)$ and therefore $\widetilde{D}$ is a $\widetilde{W}$-symmetric operator with $\deg(F_{j})\leq j$. Consequently, 
$\widetilde{D}$ belongs to $\mathcal{S}(\widetilde{W})$.
\end{proof}

\begin{remark}
    The result in Theorem \ref{algebras} may not be true if the right Fourier algebra of $\widetilde{W}$ is not diagonal. For example, if  $\widetilde{W}(x) = e^{-x^{2}}I$ and  $W(x) = T(x)e^{-x^{2}}T(x)^{\ast}$, where $T(x) = \begin{psmallmatrix} 1 & ax \\0 & 1 \end{psmallmatrix}$, then the differential operator 
    $D = \partial^{2} \begin{psmallmatrix}- 1 & ax \\ 0 & 0  \end{psmallmatrix} + \partial \begin{psmallmatrix} 0 & \frac{2}{a} \\ - \frac{2}{a} & 2x\end{psmallmatrix} + \begin{psmallmatrix} 0 & 0 \\ 0 & \frac{4}{a^{2}} \end{psmallmatrix}\in \mathcal{D}(W)$, but  it does not belong to $T(x)\mathcal{D}(e^{-x^{2}}I)T(x)^{-1}$.
\end{remark}

The following result ensures that the weight matrices $W_H$, $W_L$, and $W_J$ introduced in this paper cannot be reduced to weight matrices of smaller size. This irreducibility is a key property that further distinguishes these weights and emphasizes their significance in the context of the Matrix Bochner Problem.

\begin{prop}\label{irreducibleweights}
    The weights $W_H, W_L$, and $W_J$ are irreducible weights.
\end{prop}
\begin{proof}
    If $W$ is a reducible, then $\mathcal{D}(W)$ contains non-scalar operators of order zero (\cite{TZ18}*{Theorem 4.3}). To explore this, let $D \in \mathcal{D}(W)$ be an operator of order zero, that is  $D = F_{0} \in \operatorname{Mat}_{N}(\mathbb{C})$. By Proposition \ref{algebras}, we have that $F_{0} = T(x)G_{0}T(x)^{-1}$ for some differential operator of order zero $G_{0} \in \mathcal{D}(\widetilde{W})$. According to Theorem \ref{diag}, we obtain that $G_{0} = \operatorname{diag}(g_{1},\ldots,g_{N})$. Thus, it follows that 
     $$F_{0} = G_{0} + \sum_{j=1}^{[N/2]}a_{2j-1}x(g_{2j}-g_{2j-1})E_{2j-1,2j} + \sum_{j=1}^{[(n-1)/2]}a_{2j}x(g_{2j}-g_{2j+1})E_{2j+1,2j}.$$
     From this equality, we obtain that $g_{2j} = g_{2j-1}$ for all $1\leq j \leq [N/2]$ and $g_{2j+1} = g_{2j}$ for all $1 \leq j \leq [(N-1)/2]$. Consequently, we have that $F_{0} = G_{0} = g_{1}I$. Therefore, $\mathcal{D}(W)$ does not contain non-scalar differential operators of order zero. 
\end{proof}

Now, we will study the algebraic structure of the algebra $\mathcal{D}(W)$ associated with the weights $W = W_H, W_L,$ and $W_J$. We will prove that the leading coefficient of any differential operator in $\mathcal{D}(W)$ is a multiple of a power of $\rho(x)I$, where $\rho(x)$ takes the values $1$, $x$, or $1-x^{2}$ depending on whether $W$ corresponds to the Hermite, Laguerre, or Jacobi weights, respectively. This fact plays a crucial role in proving Theorem \ref{W-noDarboux}, which establishes that the weights $W = W_H, W_L,$ and $W_J$ cannot be 
  Darboux transformations of direct sums of scalar weights. Later, in Section \ref{sect solutions}, we will prove that the algebra $\mathcal{D}(W)$ is a polynomial algebra in a (second-order) differential operator $D$.

\begin{prop}\label{weight}
    If  $W=W_H, W_L$ or $W_J$ then we have that $\mathcal D(W)$ is a commutative algebra and any differential operator in $\mathcal{D}(W)$ is of even order.
\end{prop}
\begin{proof} 
    Both statements follow easily from Proposition \ref{algebras}, by using that $\mathcal{D}(\widetilde{W})$ is a commutative algebra and that any differential operator in $\mathcal{D}(\widetilde{W})$ is of even order.
 \end{proof}
 
Now, our goal is to determine the leading coefficient of a differential operator in the algebra $\mathcal{D}(W)$. We will prove that 
it is a multiple of $\rho(x)^m$, where $\rho$ is either $1$, $x$, or $1-x^{2}$ depending on whether the weight corresponds to Hermite, Laguerre, or Jacobi type, respectively.

The following technical result will be useful for the proof of Theorem \ref{leading}.

\begin{lemma}\label{tec}
    Let $\widetilde W= \widetilde W_H$ or $\widetilde W_L$.  Let $\mathcal{A} = \sum_{j=0}^{2m}\partial^{j}F_{j} $ be a differential operator of order $2m$ in the algebra $ \mathcal{D}(\widetilde{W})$. Then the coefficients of the differential operator $\mathcal A$ satisfy  
    $$F_{2m}(x)=  \operatorname{diag} (c_{1}\rho(x)^{m},\ldots,c_{N}\rho(x)^{m}),\quad \text{  for some } \, c_{1},\ldots,c_{N} \in \mathbb{C},  $$
    $$F_{m}(x) = \operatorname{diag}(c_{1}(-2)^{m}x^{m} + q_{1}(x), \ldots, c_{N}(-2)^{m}x^{m}+q_{N}(x)), \qquad \text {for } \widetilde W= \widetilde W_H,$$
   $$F_{m}(x) = \operatorname{diag}(c_{1}(-1)^{m}x^{m} + q_{1}(x), \ldots, c_{N}(-1)^{m}x^{m}+q_{N}(x)), \qquad \text {for } \widetilde W= \widetilde W_L,$$
    for some    polynomials $q_{i}$ of degree less than $m$,
    and $$\deg(F_{j})<j \, , \quad \text{ for  } \quad  m+1 \leq j \leq 2m. $$ 
\end{lemma}

\begin{proof} 

By Corollary \ref{alg diag}, we have that $\mathcal{A} = \operatorname{diag}(r_{1}(\delta_{1}),\cdots,r_{N}(\delta_{N}))$, with $r_{1},\ldots,r_{N} \in \mathbb{C}[x]$, and $\delta_{i} = \partial^{2}\rho(x) + \partial (sx+t)$, where $\rho(x) = 1$, $s = -2$, and $t = -2b_{i}$ for the Hermite-type weight, and $\rho(x) = x$, $s = -1$, and $t = \alpha_{i} + 1$ for the Laguerre-type weight. It is sufficient to prove the statement for $N=1$, and $r_{1}(x) = x^{m}$. 
 We recall that given  a differential operator $\nu =\sum_{j=0}^{m}\partial^{j}f_{j} \in \mathcal{D}(w)$, where $f_{j}(x) = \sum_{k=0}^{j}c_{j,k}x^{k}$, we have that the sequence of eigenvalues of $\nu$ satisfies $\Lambda_{n}(\nu) = \sum_{j=0}^{m}c_{j,j}[n]_{j}$, where $[n]_{j} = n(n-1)\cdots(n-j+1)$. In particular $\Lambda_n(\delta_{1})= sn$ because $\deg(\rho(x))\leq 1$. 

By Proposition \ref{prop2.8-GT}, we have that 
    \begin{align*}
\Lambda \colon \mathcal{D}(w) &\longrightarrow \Lambda(\mathcal{D}(w)) \\
\nu &\mapsto \{\Lambda_{n}(\nu)\}_{n}
\end{align*}
is an algebra isomorphism. Hence $\Lambda_{m}(\delta_{1}^{m})= \Lambda_{n}(\delta_{1})^{m} = ( s n)^{m}$.
We write $\delta_{1}^{m} = \sum_{j=0}^{2m}\partial^{j}f_{j}$ with $f_{j}(x) = \sum_{k=0}^{j}c_{j,k}x^{k}$, then
$$(sn)^m= \Lambda_{m}(\delta_{1}^{m}) = \sum_{j=0}^{2m}c_{j,j}[n]_{j}.$$ 
Now, observe that the left-hand side of the above equation is polynomial in $n$ of degree $m$, therefore 
$c_{j,j} = 0$ for all $m < j \leq 2m$, and  $c_{m,m} = s^{m}$. This completes the proof.
\end{proof}

We have proved that each differential operator in $\mathcal{D}(W)$ is of even order. In the following theorem, we characterize the leading coefficient of these operators. 

\begin{thm}\label{leading} 
Let $W$ be one of the weights $W_H$, $W_L$, or $W_J$ and let $\rho(x)$ be equal to  $1$, $x$, or $1-x^2$, corresponding to the Hermite, Laguerre, or Jacobi weight, respectively.
Then, if $\mathcal{A} \in \mathcal{D}(W)$ is a differential operator of order $2m$, the leading coefficient of $\mathcal{A}$ is a multiple of $\rho(x)^m$. 
\end{thm}

\begin{proof}

If $\mathcal{A} = \sum_{j=0}^{2m}\partial^{j}G_{j} \in \mathcal{D}(W)$ is a differential operator of order $2m$, then according to Proposition \ref{algebras}, there exists a differential operator $\mathcal{B} = \sum_{j=0}^{2m}\partial^{j}F_{j}\in \mathcal{D}(\widetilde{W})$ of order $2m$ such that $\mathcal{A} = T(x)\mathcal{B}T(x)^{-1}$. 

We proceed with the proof by considering two separate cases: Jacobi-type on one hand, and Hermite-type and Laguerre-type on the other hand.

\smallskip
    {\em Jacobi-type:} 
From Remark \ref{leadingDWtilde}, it follows that $F_{2m}(x) = \operatorname{diag}(c_{1}(1-x^{2})^{m},\ldots,c_{N}(1-x^{2})^{m})$ for some $c_{1},\ldots,c_{N} \in \mathbb{C}$.  By \eqref{coeff}, we have 
$$G_{2m}(x) = T(x)F_{2m}T(x)^{-1} = F_{2m}(x) + (AF_{2m}(x) - F_{2m}(x)A)x.$$ 
Since $\deg(G_{2m})$ and $\deg(F_{2m})$ are less than of equal to $2m$, then $(AF_{2m}-F_{2m}A)x$ must be of degree less than or equal to $2m$. We have
\begin{equation*}
     \begin{split}
        x [A,F_{2m}(x)] & =x (1-x^{2})^{m}  \Big( \sum_{j=1}^{[N/2]}{a_{2j-1}} (c_{2j}-c_{2j-1}) E_{2j-1,2j} + \sum_{j=1}^{[(N-1)/2]}{a_{2j}}
         (c_{2j}-c_{2j+1})) E_{2j+1,2j} \Big) .
     \end{split}
 \end{equation*}
From here, we obtain that $c_{2j}-c_{2j-1} = 0$ and $c_{2j}-c_{2j+1} = 0$ for all $j$. Therefore, $c_{j} = c_{i}$ for all $1\leq i,j \leq N$, and $F_{2m}(x) = c_{1}(1-x^{2})^{m}I = G_{2m}.$

\smallskip 
 {\em Hermite-type and Laguerre-type:}  
 We have that $$F_{2m}(x) = \operatorname{diag}(c_{1}\rho(x)^{m},\ldots,c_{N}\rho(x)^{m}), \quad \text{ for some $c_{1},\ldots,c_{N}\in \mathbb{C}$,}$$ with $\rho(x) = 1$ for Hermite-type and $\rho(x) = x$ for Laguerre-type. 
 
 From Lemma \ref{tec} we know  that $\deg(F_{j}) < j$ for all $m < j \leq 2m$, and $$F_{m}(x) = \operatorname{diag}\big(c_{1}s^{m}x^{m} + q_{1}(x),\ldots,c_{N}s^{m}x^{m}+q_{N}(x)\big)$$ with $q_{j}$ a polynomial, $\deg(q_j)<m$, $s=-2$ for Hermite-type and $s=-1$ for Laguerre-type. 
 
 By \eqref{coeff}, it follows that $G_{m} = F_{m} + [A,F_{m}]x + (m+1)AF_{m+1}$. Since $\deg(G_{m})$, $\deg(F_{m})$ and $\deg(AF_{m+1})$ are less than or equal to $m$, we obtain that $[A,F_{m}]$ must be of degree less than or equal to $m-1$. We have that
 \begin{equation*}
     \begin{split}
         [A,F_{m}(x)] & = \sum_{j=1}^{[N/2]}{a_{2j-1}}\big (s^{m}x^{m}(c_{2j}-c_{2j-1}) + (q_{2j}-q_{2j-1}) \big) E_{2j-1,2j} \\
         & \quad + \sum_{j=1}^{[(N-1)/2]}{a_{2j}}\big(s^{m}x^{m}(c_{2j}-c_{2j+1})+(q_{2j}-q_{2j+1}) \big)E_{2j+1,2j}. 
     \end{split}
 \end{equation*}
 Thus  $c_{2j}-c_{2j-1} = 0$ and $c_{2j}-c_{2j+1} = 0$ for all $j$, i.e. $c_{j} = c_{i}$ for all $1\leq i,j \leq N$. Therefore $F_{2m}(x) = c_{1}\rho(x)^{m}I$ and  $G_{2m}(x) = T(x)F_{2m}(x)T^{-1}(x) = c_{1}\rho(x)^{m}I$.
\end{proof}

\

We are now in a position to provide a proof of the following theorem, which was stated in the Introduction.

\begin{thm}\label{W-noDarboux}
 The irreducible weight matrices  $W_H, W_L$, and $W_J$ are not Darboux transformations of any direct sum of scalar weights.
\end{thm}
\begin{proof}
Let $W = W_{H}, \, W_{L},$ or $W_{J}$, and $\widetilde{W}$ be a direct sum of scalar weights. Let $P_{n}(x)$ and $\widetilde{P}_{n}(x)$ be the sequences of monic orthogonal polynomials for $W$ and $\widetilde{W}$ respectively.
If $W$ is a bispectral Darboux transformation of $\widetilde{W}$, then there exist differential operators $\mathfrak{D}, \widetilde{\mathfrak{D}} \in \mathcal{F}_{R}(W)$, polynomials $F,\widetilde{F}$, and sequences of matrices $C_n, \widetilde{C}_n$ which are nonsingular for almost every $n$ such that 
$$C_n\widetilde{P}_n(x) = P_n(x) \cdot \mathfrak{D}F(x)^{-1} \text{ and } \quad  \widetilde{C}_n P_n(x) = \widetilde{P}_n(x) \cdot \widetilde{F}(x)^{-1}\widetilde{\mathfrak{D}}.$$
The zero-order differential operator $E_{11}$ belongs to $\mathcal{D}(\widetilde{W})$ and we have that 
$$ P_n(x)\cdot \mathfrak{D}F(x)^{-1}E_{11}\widetilde{F}(x)^{-1}\widetilde{\mathfrak{D}}= C_nE_{11}\widetilde{C}_nP_n(x).$$
Therefore, the nonzero differential operator $\mathfrak{D}F(x)^{-1}E_{11}\widetilde{F}(x)^{-1}\widetilde{\mathfrak{D}}$ belongs to $\mathcal{D}(W)$, but its leading coefficient is a singular matrix, which contradicts Theorem \ref{leading}.
\end{proof}



\section{Solutions of the Matrix Bochner Problem} \label{sect solutions}

In this section, we will prove that the weights \(W = W_H, W_L,\) and \(W_J\) are solutions to the Matrix Bochner Problem, meaning they admit a second-order differential operator in the algebra \(\mathcal{D}(W)\). Furthermore, we will show that for these weights \(W\), the algebra \(\mathcal{D}(W)\) is always a polynomial algebra in a second-order differential operator.

Recall that we are considering weights of the form $W(x)=T(x)\widetilde W(x) T(x)^*$, where $\widetilde{W}(x) = w_{1}(x) \oplus \cdots \oplus w_{N}(x)$ is a direct sum of classical scalar weights of the same type. From Proposition \ref{fourier} and  Theorem \ref{algebras} we have the following relationship between the algebras 
 $$\mathcal{F}_R(W) = T(x)\mathcal{F}_R(\widetilde{W})T^{-1}(x), \quad \text{ and } \quad \mathcal{D}(W) \subset T(x)\mathcal{D}(\widetilde{W})T^{-1}(x).$$

\smallskip
The classical Hermite, Laguerre, or Jacobi operators $\delta$ can be written  as 
 \begin{equation*}
     \delta= \partial ^2 \rho(x) + \partial \, \tfrac{ (w(x) \rho(x))' }{w(x)},
 \end{equation*} 
where $\rho(x)=1, \, x, \text{ or } (1-x^2)$, depending on whether $w(x)$ is the classical Hermite, Laguerre, or Jacobi weight, respectively. Therefore,  the diagonal differential operator 
\begin{equation}\label{Dtilde}
\widetilde D= \partial^{2}\rho(x) + \partial \, \left (\widetilde W(x)\rho(x) \right )'\, \widetilde W(x)^{-1} 
\end{equation}
is a $\widetilde W$-symmetric operator in the algebra $\mathcal D(\widetilde W)$. 
Unfortunately, the operator $T(x)\,\widetilde D\, T(x)^{-1}$  does not belong to the algebra  $\mathcal D(W)$, even when it 
 is a $W$-symmetric operator in the Fourier algebra $\mathcal F_R(W)$.
 However, we can find another second-order differential operator in the algebra $\mathcal D(W)$, closely related with $T(x)\,\widetilde D\, T(x)^{-1}$. Indeed, we have the following result.

  \begin{prop}  \label{solutionDgeneral}
      Let $ W =  W_{H}, W_{L}$, or $  W_{J}$ be the weight matrices introduced in Definition \ref{W-definition}. Then there exists a constant matrix $\widetilde K$ such that 
         $$ D = T(x)\left(\partial^{2}\rho(x) + \partial (\widetilde W(x)\rho(x))'\widetilde W(x)^{-1} + \widetilde K\right) T(x)^{-1} \in \mathcal D(W). $$   
Explicitly, if $K= \sum_{j=1}^{[N/2]} E_{2j,2j}$ then $\widetilde K= 2K, \, K$, or $(\alpha_1+\beta_1)K$ according to  whether $W=W_{H}, W_{L}$, or $W_{J}$, respectively.  
  \end{prop}

\begin{proof} From Table \ref{weights-table}, we have that 
diagonal weight $\widetilde W$ satisfies 
$(\widetilde W(x)\rho(x))'\widetilde W(x)^{-1} = B-xL$  
where  the diagonal matrices $B$ and $L$ are given by 
\begin{align*}
   B & =\textstyle \sum_{j=1}^N 2b_j E_{j,j} , \qquad  \qquad  
   L = 2\operatorname{Id},      &\text{ for  Hermite weights.}  \\ 
  B &=\textstyle\sum_{j=1}^N (\alpha_j+1) E_{j,j} , \qquad 
  L= \operatorname{Id},  &  \text{ for Laguerre weights.} \\
    B &=\textstyle\sum_{j=1}^N (\beta_j-\alpha_j) E_{j,j} , \quad \; \; 
    L= 2\sum_{j=1}^N (\alpha_j+\beta_j+2 ) E_{j,j},   & \text{for  Jacobi weights.} 
\end{align*}

Let $\widetilde{D}$ be the $\widetilde W$-symmetric differential operator given  in \eqref{Dtilde}, and $T(x) = e^{Ax} = I + Ax$. Then the differential operator $D = T(x)\left(\widetilde{D} + \widetilde{K}\right)T(x)^{-1}$ is $W$-symmetric, and it is of the form 
$D=  \partial^{2}\rho(x) + \partial F_1(x) + F_0(x),$ with 
\begin{align*}
    F_1(x)& =  2\rho(x) A +B-xL+x[A,B-xL]  \quad \text { and } \quad F_0(x)=  A(B-xL)+\widetilde K+ x[A,\widetilde K].
\end{align*}
The operator $D$ belongs to the algebra $\mathcal D(W)$ if and only if  $F_1(x)$ and $F_0(x)$  are polynomials of degree at most 1 and 0 respectively.

For  Hermite and Laguerre weights, it is clear that $F_1(x)$ is always a polynomial of degree one, and $F_0$ is a constant matrix since $AL= [A,\widetilde K]$.  

In the case of  Jacobi weights, the parameters $\alpha_j, \, \beta_j$ satisfy the extra conditions $\alpha_j+ \beta_j +1+(-1)^j =\alpha_1+\beta_1$, for all $ j$. Hence $L=(\alpha_1+\beta_1+2)I- 2K $. In particular, $[A,L]=-2[A,J]= -2A$ and $AL=  [A, \widetilde J]$. Therefore, 
$$F_1(x)= 2A+B + x\big( [A,B]+2 J -(\alpha_1+\beta_1+2)I \big) \quad \text{ and }
\quad F_0(x)= AB+(\alpha_1+\beta_1) J. $$
Thus, the operator $D$ belongs to the algebra $\mathcal D(W)$. 
\end{proof}

By making explicit the expressions of the coefficients of the differential operator $D$ in Proposition \ref{solutionDgeneral}, corresponding to each family we obtain the following result.

\begin{thm}\label{boch sol} 
Let   $W=W_H, W_L$ or $W_J$ be the weight matrices introduced in Definition \ref{W-definition} and let  
$$ K=\sum_{j=1}^{[N/2]} E_{2j,2j} , \quad B_H= \textstyle \sum_{j=1}^N 2b_j E_{j,j} \, , \quad  B_L=\textstyle \sum_{j=1}^{N} (\alpha_j+1) E_{j,j}\, , 
\quad   B_J=\textstyle \sum_{j=1}^{N} (\beta_{j}-\alpha_j) E_{j,j}.$$
Then 
 \begin{align*}
  D_H &= \partial ^2 I+ \partial \Big( 2A+B_H-x+x[A,B_H]\Big) + AB_H+ 2K \in \mathcal{D}(W_{H}),\\ 
    D_L &= \partial ^2 x I+ \partial \Big( B_L-x+2xA+x[A,B_L]\Big) + AB_L+ K\in \mathcal{D}(W_{L}),\\
        D_J & = \partial ^2 (1-x^2) I+ \partial \Big(2A+B_J + x\big( [A,B_J]+2 K -(\alpha_1+\beta_1+2)I\big)\Big) + AB_J+(\alpha_1+\beta_1)K 
    \end{align*}
     belongs to $\mathcal{D}(W_{J})$.

\end{thm}

Upon establishing the existence of a second-order differential operator in the algebra and in light of our knowledge of the leading coefficients of each operator in Theorem \ref{leading}, we obtain the structure of the algebra $\mathcal D(W)$.

\begin{thm} \label{poly alg}
    Let $W$ be one of the weights $W_{H}$, $W_{L}$, or $W_{J}$ and let $D$ be  the  second-order 
    differential operator  $D_H, D_L$ or $D_J$ in $\mathcal{D}(W)$ respectively.
    Then $\mathcal{D}(W)$ is a polynomial algebra in $D$. 
\end{thm}
\begin{proof} 
Let $\mathcal{A} \in \mathcal{D}(W)$ a differential operator of order $2m$. 
 We proceed by induction on $m$. For $m = 0$, the statement follows from Proposition \ref{irreducibleweights}.
 For $m > 0$, we write $\mathcal{A} = \sum_{j=0}^{2m}\partial^{j}F_{j}$.  From Proposition \ref{leading} we have that $F_{2m}(x) = c\rho(x)^{m}$ for some $c \in \mathbb{C}$. The leading coefficient of the second-order differential operator $D$ in $\mathcal D(W)$ is equal to $\rho(x)$. Hence, the differential operator $\mathcal{A} - cD^{m}$ belongs to $\mathcal{D}(W)$ and has order less than $2m$. This completes the proof.
\end{proof}

\section{Explicit examples of size  $2\times 2 $} 

For the benefit of the reader, 
we give explicitly  the weights $W_{H}$, $W_{L}$ and $W_{J}$ for $N=2$, along with their corresponding second-order differential operators in the algebra $\mathcal D(W)$.

\subsection{Hermite weights of size 2}
     For $N = 2$, $a,b \not=0$, we have 
     $$W_{H}(x) = \begin{pmatrix} 1 & ax \\ 0 & 1 \end{pmatrix} \begin{pmatrix} e^{-x^{2}+2bx} & 0 \\ 0 & e^{-x^{2}} \end{pmatrix} \begin{pmatrix} 1 & 0 \\ ax & 1 \end{pmatrix} = e^{-x^{2}} \begin{pmatrix} e^{2bx} + a^{2}x^{2} && ax \\ ax && 1 \end{pmatrix}$$
     together with the second-order differential operator
     $$D = \partial^{2} I + \partial \begin{pmatrix}-2x + 2b && -2abx + 2a \\ 0 && -2x  \end{pmatrix} + \begin{pmatrix} 0 && 0 \\ 0 && 2 \end{pmatrix} \in \mathcal{D}(W_{H}).$$

\subsection{Laguerre weights of size 2}
    For $N = 2$, $\alpha,\beta > -1$, $a \not = 0$, $\beta - \alpha \notin \mathbb{Z}$ we have 
    $$W_{L}(x) = \begin{pmatrix} 1 & ax \\ 0 & 1 \end{pmatrix} \begin{pmatrix} e^{-x}x^{\alpha} & 0 \\ 0 & e^{-x}x^{\beta} \end{pmatrix} \begin{pmatrix} 1 & 0 \\ ax & 1 \end{pmatrix} =e^{-x} x^{\beta} \begin{pmatrix} x^{\alpha-\beta} + a^{2}x^{2} && ax \\ ax && 1 \end{pmatrix}$$
    together with the second-order differential operator
    $$D = \partial^{2} xI + \partial \begin{pmatrix}\alpha + 1 -x&& ax(2+\beta-\alpha) \\ 0 && \beta + 1 -x  \end{pmatrix} + \begin{pmatrix} 0 && a(\beta+1) \\ 0 && 1 \end{pmatrix}\in \mathcal{D}(W_{L}).$$

\subsection{Jacobi weights of size 2}
For $\alpha_{1}, \, \alpha_{2}, \, \beta_{1}, \, \beta_{2} > -1$, $a\not=0$, and $\alpha_1-\alpha_2\not \in \mathbb Z$ or $\beta_1-\beta_2\not \in \mathbb Z$,
we have the irreducible weight 

\begin{equation*}
    \begin{split}
        W_{J} & = \begin{pmatrix} 1 & ax \\ 0 & 1 \end{pmatrix} \begin{pmatrix} (1-x)^{\alpha_{1}}(1+x)^{\beta_{1}} & 0 \\ 0 & (1-x)^{\alpha_{2}}(1+x)^{\beta_{2}} \end{pmatrix} \begin{pmatrix} 1 & 0 \\ ax & 1 \end{pmatrix} \\
        & =(1-x)^{\alpha_2}(1+x)^{\beta_2}\begin{pmatrix}
    (1-x)^{\alpha_1-\alpha_2}(1+x)^{\beta_1-\beta_2} +a^2 x^2 & & ax\\ ax & & 1
\end{pmatrix}.
    \end{split}
\end{equation*}

If $\alpha_1+\beta_1=\alpha_2+\beta_2+2$, then the following differential  operator $D$ belongs to the algebra $\mathcal D(W_{J})$

\begin{align*}
D&= \partial^{2} (1-x^{2})I + \partial  \begin{pmatrix} \beta_{1} - \alpha_{1} -x(\alpha_{1} + \beta_{1} + 2) && a(2 + x (\beta_{2} - \beta_{1} + \alpha_{1} - \alpha_{2})) \\ 0 && \beta_{2} - \alpha_{2} - x(\alpha_{2} + \beta_{2} + 2)  \end{pmatrix} \\ 
& \qquad + \begin{pmatrix} 0 && a(\beta_{2} - \alpha_{2}) \\ 0 && \alpha_{1}+\beta_{1}  \end{pmatrix}. 
\end{align*}

\end{document}